\documentclass[12pt]{amsart}
\usepackage{amsmath}
\usepackage{tikz}
\usetikzlibrary{calc} 
\usepackage{amsthm}
\usepackage{amssymb}
\usepackage{graphicx}
\usepackage{blindtext}
\usepackage{hyperref}
\numberwithin{equation}{section}

\usepackage[margin=2.75cm]{geometry}
\makeatletter
\@namedef{subjclassname@2020}{\textup{2020} Mathematics Subject Classification}
\makeatother

\newtheorem{theorem}{Theorem}[section]

\newtheorem{lemma}[theorem]{Lemma}

\newtheorem{corollary}[theorem]{Corollary}
\newtheorem{definition}[theorem]{Definition}

\newtheorem{remark}[theorem]{Remark}
\date{\today}
\begin{document}


\title[Asymptotic Behavior of the Non-resonance Eigenvalues]{Asymptotic Behavior of the Non-resonance Eigenvalues of the Fractional Schrödinger Operator with Neumann Condition}


\author[S. Karakılıç]{Sedef Karakılıç} 
\author[S.~Özcan]{Sedef Özcan} 
\address{Sedef Karakılıç, Dokuz Eylül University, Faculty of Science, Department of Mathematics, Izmir, Turkey}
\address{Sedef Özcan, Dokuz Eylül University, Faculty of Science, Department of Mathematics, Izmir, Turkey}



\maketitle

\begin{abstract}
We present an analytical investigation of the asymptotic behavior of non-resonance eigenvalues for the fractional Schrödinger operator under homogeneous Neumann boundary conditions. Our findings reveal an intriguing convergence: as the system evolves, the eigenvalues of the fractional Schrödinger operator increasingly resemble those of the fractional Laplace operator. By deriving a precise asymptotic formula, we provide new insights into the spectral properties of these operators, highlighting their deeper connections and potential applications in mathematical physics.


\end{abstract} 



\section{Introduction}
\setcounter{section}{1} \setcounter{equation}{0}
The asymptotic behavior of eigenvalues of self-adjoint operators is fundamental in quantum mechanics, spectral geometry, and mathematical physics, providing insights into energy distribution, wave propagation, and geometric properties of manifolds \cite{29}. It plays a crucial role in understanding the semiclassical limit, Weyl’s law, and quantum chaos, with applications in various physical models, including random Schrödinger operators and wave localization phenomena \cite{23}.
In particular, for the Schrödinger operator, eigenvalue asymptotics describe the density of states, spectral properties of quantum systems, and localization effects in disordered media \cite{30}. 

The \emph{fractional Neumann Laplace operator} has significant applications in areas where spectral properties are essential. In \emph{numerical analysis}, efficient spectral methods have been developed to approximate solutions involving the spectral fractional Laplacian, particularly on rectangular domains. These methods utilize Fourier-like basis functions to discretize the operator, facilitating the resolution of partial differential equations with spectral fractional Laplacians  on rectangular domains  \cite{27}.
In \emph{turbulence modeling}, the spectral fractional Laplacian offers a framework for analyzing anomalous diffusion and energy transport across scales in complex systems, such as dusty plasma monolayers. This approach aids in understanding the spectral characteristics of turbulence in such media \cite{28}.

The fractional Schrödinger operators are particularly important in anomalous diffusion and nonlocal quantum interactions, with applications in topological insulators, quantum Hall effects, and fractional quantum mechanics \cite{32,31}.

Let $K=[0,a_{1}]\times[0,a_{2}]\times\cdots\times[0,a_{d}]$
and  $-\Delta_{\mathsf{NEU}}$   denote the Laplace
operator  
\begin{equation}\label{NeumLap}
\Delta_{\mathsf{NEU}}  =\frac{\partial^2}{\partial x_1^2}+\cdots+\frac{\partial^2}{\partial x_d^2}
\end{equation}
defined in $L_2(K)$ with the Neumann boundary conditions
\begin{equation}\label{NeumBC} 
\frac{\partial u}{\partial n} \Big|_{\partial K}=0,\end{equation} 
where $x=(x_{1},x_{2},\cdots,x_{d})\in
R^{d} $, $ d\geq2 $, and
$\partial K $ is the boundary of the domain $ K $,
$\frac{\partial}{\partial n}$ is the differentiation along the
outward normal $n$ of $\partial K.$ 
	In this paper, we examine the fractional Schrödinger operator
\begin{equation}\label{shrod}
\mathsf{H}_{\mathsf{NEU}}(\ell,q)u=(-\Delta_{\mathsf{NEU}})^{\ell} u+q(x)u,
\end{equation}
where  $q(x)\in L_{2}(K)$ and  $\frac{1}{2}< \ell < 1$.
	Due to its physical significance, the most substantial advancements have been made for the Schrödinger operator, specifically when \( \ell = 1 \).
    In this context, O. A. Veliev made a significant contribution to the study of the periodic Schrödinger operator (concerning an arbitrary lattice and arbitrary dimension $d\geq 2$) in papers \cite{8,9,10,11} by being the first to classify the large Bloch eigenvalues of the Laplace operator into two distinct groups: non-resonance and resonance eigenvalues. His profound insights led to the derivation of asymptotic formulas that elegantly describe the perturbations within each group, which played a crucial role in affirming the Bethe-Sommerfeld conjecture. His book \cite{12} is one of the key contemporary references on the multidimensional periodic Schrödinger operator, where his method, previous results of \cite{8,9,10,11}, and more are examined in detail. Additional proofs of these asymptotic formulas for quasiperiodic boundary conditions in two- and three-dimensional cases can be found in \cite{2,3,6,7}. For the Schrödinger operator with periodic boundary conditions, asymptotic formulas for the eigenvalues were obtained in \cite{4}. In papers \cite{18,13,1}, it is derived that the eigenvalue formulas for the Schrödinger operator with Dirichlet and Neumann boundary conditions on a \(d\)-dimensional parallelepiped for any \(d \geq 2\).

 The aim of this paper is to examine the behavior of the eigenvalues of the fractional Schrödinger operator with Neumann boundary conditions as they grow larger, adapting O. A. Veliev’s method by relying on the idea presented in \cite[Remark 1]{9}. Particularly;  for arbitrary \( \frac{1}{2}<\ell<1\) and any \(d \geq 2\), the primary objective of this paper is to achieve that for non-resonance eigenvalues  $|\beta|^{2\ell}$ of the $(-\Delta_{\mathsf{NEU}})^{\ell}$ where $|\beta|>>1$,
 	there exists an eigenvalue $\xi_{N}$ of $\mathsf{H}_{\mathsf{NEU}}$ satisfying
 	\begin{equation}
 	\xi_{N}=|\beta|^{2\ell} + F_{k-1} +O(r^{-k\alpha_1(\ell)}), \hspace{.3cm} \forall k=1,2,\cdots,p-c,
 \end{equation}
 where $F_k$'s depend on the Fourier coefficients of the potential $q(x)$ satisfying the assumption \eqref{cond2} with $|F_k|=O(r(\ell))$ (see Section \ref{asymp}). 
	
	Let us outline the structure of this article. In Section \ref{sec2}, we will present essential concepts. In Section \ref{sec3}, we derive iteration series, and owing to it, we accomplish to prove that the non-resonance eigenvalues of the fractional Schrödinger operator and the fractional Laplace operator are asymptotically close to each other. 
\section{Preliminaries}\label{sec2}
\setcounter{section}{2} \setcounter{equation}{0}
	Let  \(K \equiv [0, a_1] \times [0, a_2] \times \cdots \times [0, a_d]\) be a rectangular subset of $\mathbb{R}^{d}$, where $d\geq 2$ and \(\partial K\) denote its boundary. We denote by \(\frac{\partial}{\partial n}\) the differentiation along the outward normal \(n\) of \(\partial K\). 

We will take the domain of $-\Delta_{\mathsf{NEU}}$ defined by \eqref{NeumLap} and \eqref{NeumBC} to be the subspace of $L_{2}(K)$ consisting of smooth functions with Neumann boundary values, and denote it by $\mathcal{D}(-\Delta_{\mathsf{NEU}})$. The eigenvalues of $-\Delta_{\mathsf{NEU}}$ are
	$|\beta|^{2}$ with the corresponding eigenfunctions
		\begin{equation*}
		v_{\beta}(x)=\cos \beta^{1}x_1 \cdots \cos \beta^{d}x_{d}, 
	\end{equation*}
where
\begin{align*}
	\beta \in
\mathcal{B}^{+}=&\left\lbrace \beta=(\beta^1,\beta^2,\cdots,\beta^d) : \beta^i=\frac{n_i \pi}{a_i}, n_i \in \mathbb{Z}^{+0}, i=1,2,\cdots,d\right\rbrace .
\end{align*}
The system $ \left\lbrace v_{\beta} \right\rbrace_{\beta \in\mathcal{B}^{+}} $ is known to be an orthogonal basis for $L_{2}(K)$. Therefore, any $u(x) \in L_{2}(K)$ has the Fourier series expansion
		\begin{align}\label{qx1}
				u(x)= \sum_{\beta \in\mathcal{B}^{+}}u_{\beta}v_{\beta}(x)
		\end{align}
		where $u_{\beta}=\frac{1}{\mu(K)}(q, v_{\beta})$.

The spectral characterization of the fractional Neumann Laplacian $(-\Delta_{\mathsf{NEU}})^{\ell}$, which is based on the spectral theorem (see \cite[Theorem VIII.6]{23}), starts with $-\Delta_{\mathsf{NEU}}$ and is followed by taking its spectral power. The spectral fractional Neumann Laplacian $(-\Delta_{\mathsf{NEU}})^{\ell}$ for $\frac{1}{2}<\ell<1$ is defined in \cite[Section 2.5.1]{17} as
\[
    (-\Delta_{\mathsf{NEU}})^{\ell} u = \sum_{\beta \in\mathcal{B}^{+}}|{\beta}|^{2\ell}u_{\beta}v_{\beta}(x),\; \; \; u_{\beta}=\frac{1}{\mu (K)}(u,v_{\beta}),
\]
for every $u \in \mathcal{D}(-\Delta_{\mathsf{NEU}})$. Equivalently, via the spectral mapping theorem (see also \cite[Proposition 10.3]{19})  
	 the eigenvalues of $(-\Delta_{\mathsf{NEU}})^{\ell}$ are
	$|\beta|^{2\ell}$ with the corresponding eigenfunctions
	$v_{\beta}(x)$, $\beta \in\mathcal{B}^{+}$. 
	Moreover, $(-\Delta_{\mathsf{NEU}})^{\ell}$ preserves the self-adjointness and positivity properties of \( \Delta_N \). 

	The fractional Schrödinger operator $\mathsf{H}_{\mathsf{NEU}}(\ell,q)$ is defined as
		\begin{equation}\label{perturbedoperator}
		\mathsf{H}_{\mathsf{NEU}}(\ell,q)u(x)={(-\Delta_{\mathsf{NEU}})}^{\ell} u(x)+q(x)u(x),
		\end{equation}
		for all $u(x) \in \mathcal{D}(-\Delta_{\mathsf{NEU}})$, where $q(x)$ is a real valued function in $L_2(K)$ satisfying the condition \eqref{cond2}.
	The eigenvalues of $\mathsf{H}_{\mathsf{NEU}}(\ell,q)$ are denoted by $\xi_N$ with the corresponding normalized eigenfunctions $\chi_{N}$. We denote by $(.,.)$ and $ \left\langle .,. \right\rangle  $ the inner products in $L_{2}(K)$ and $\mathbb{R}^d$, respectively.

 By the fact \eqref{qx1} , the potential $q(x) \in L_{2}(K)$ has the Fourier series expansion
		\begin{align}\label{qx}
				q(x)= \sum_{\beta \in\mathcal{B}^{+}}q_{\beta}v_{\beta}(x)
		\end{align}
		where $q_{\beta}=\frac{1}{\mu(K)}(q, v_{\beta})$.
	Let
	\begin{align*}
\mathcal{B}=&\left\lbrace \beta=(\beta^1,\beta^2,\cdots,\beta^d) : \beta^i=\frac{n_i \pi}{a_i}, n_i \in \mathbb{Z}, i=1,2,\cdots,d\right\rbrace,\\
	A_{\beta}=&\left\lbrace \alpha
= (\alpha^{1}, \alpha^{2},\cdots,\alpha^{d}) \in\mathcal{B} :
|\alpha^{i}|=|\beta^{i}|, i=1,2,\cdots,d \right\rbrace.
	\end{align*}
		For $q(x)\in L_{2}(K)$, we have (see \cite{18})
		\begin{align}\label{qx}
			q(x)= \sum_{\beta \in\mathcal{B}}q_{\beta}v_{\beta}(x),
		\end{align}
		where $q_{\beta}=\frac{1}{|A_{\beta}|\mu(K)}(q(x),v_{\beta}(x))$ and $|A_{\beta}|$ denotes the number of vectors in $A_{\beta}$.
	Without loss of generality, we assume that $q_0=0.$
	We also assume that the potential $q(x)$ satisfies the following condition:
	\begin{equation}\label{cond2}
		\sum_{\beta\in
			\mathcal{B}}|q_{\beta}|^2(1+|\beta|^{2m})<\infty,
	\end{equation}
	where $m>\frac{(4d-1)}{2}(d+20)3^{d+1}+\frac{d}{4}3^{d}+d+1 $, which is a similar condition to in \cite{12,16}. In other words, the potential function $q(x)$ is differentiable enough (of order $m$). By \eqref{cond2}, we also have 
	\begin{equation}\label{qxorder}
		q(x)=\sum_{\beta\in \mathcal{B}(r^{\alpha(\ell)})}q_{\beta}v_{\beta}(x)+O(r^{-p\alpha(\ell)}),
	\end{equation}
	where $p=m-d, \alpha(\ell)$ is a power smaller than $\frac{1}{d+20}$. By \eqref{cond2}, we have
		\begin{equation}\label{Mqx}
		 M=\sum_{\beta\in\mathcal{B}} |q_{\beta}|<\infty.
	\end{equation}	
\section{Asymptotic Formula for the Eigenvalues in the Non-Resonance Domain}\label{sec3}
\setcounter{section}{3} \setcounter{equation}{0}
 We begin this section by dividing the set of eigenvalues $|\beta|^{2\ell}$ for $|\beta|\sim r$ into two groups:  resonance and non-resonance ones. To achieve this, we define the corresponding resonance  and non-resonance sets which are the so-called sets in \cite{8,9,10,11,12}.  $|\beta|\sim r$ means that there exist constants $c_i\in \mathbb{R},i=1,2$,, independent of $r$, such that $c_1 r< |\beta|<c_{2}r$. We will prove an asymptotic formula for the non-resonance eigenvalues. 
Let $\alpha(\ell)=  \frac{2\ell-1}{2(d+20)3^{d+1}} $ and $\alpha_{k}(\ell)=3^{k}\alpha(\ell)$.
\begin{definition}
	For $\beta \in \mathcal{B}(pr(\ell))$, let $$V_{\beta}^{\ell}(r(\ell))\equiv\left\{x\in R^{d}:\left||
	x |
	^{2\ell}- | x+\beta |^{2\ell}\right |<r(\ell)\right\},$$
	where $r(\ell)=r^{\alpha_{1}(\ell)}$. The domains $V_{\beta}^{\ell}(r(\ell))$ are referred to resonance domains and the eigenvalue $| \gamma |^{2\ell}$ is a resonance eigenvalue if $\gamma \in V_{\beta}^{\ell}(r(\ell)).$ 
	The union of all resonance domains	$$E^{\ell}_{1}(r(\ell),p)= \bigcup_{\beta \in \mathcal{B}(p r^{\alpha(\ell)})} V^{\ell}_{\beta}(r(\ell)), $$
	is said to be resonance set and its complement
	$$ U^{\ell}(r(\ell),p)= \mathbb{R}^d \backslash E^{\ell}_{1}(r(\ell),p),$$
	is called non-resonance domain. If $\beta\in U^{\ell}(r(\ell),p)$, then the eigenvalue $| \beta|^{2\ell}$ is called a non-resonance eigenvalue. 
\end{definition} 
For notational convenience, the sets $U^{\ell}(r(\ell),p)$ and  $V_{\beta}^{\ell}(r(\ell))$ are the analogues of the resonance and non-resonance sets $U_{\beta}(r^{\alpha})$ and $V_{\beta}(r^{\alpha})$, respectively, as introduced  in, forexample, \cite{12}.
The following lemma establishes a relationship between  the non-resonance domains \(U_{\beta}^{\ell}(r(\ell))\) and \(U_{\beta}(r^{\alpha_1(\ell) - 2\ell + 2})\) and  the resonance domains \(V_{\beta}^{\ell}(r^{\alpha_k(\ell)})\) and \(V_{\beta}(r^{\alpha_k(\ell) - 2\ell + 2})\). 

The set \( V_{\beta}^{\ell}(r(\ell)) \) is defined as the collection of all vectors \( x \in \mathbb{R}^d \) satisfying the condition that the absolute difference between the \( 2\ell \)-th power of the norm of \( x \) and that of \( x + \beta \) is bounded by \( r(\ell) \).
Here, the vector \( \beta \in \mathbb{R}^d \) dictates the orientation of the set, whereas the parameter \( r(\ell) > 0 \) governs the admissible deviation. Consequently, the set can be interpreted as a region in \( \mathbb{R}^d \) that is shaped along the direction of \( \beta \), with its size determined by the parameters \( \ell \) and \( r(\ell) \).

\begin{figure}[htpb]
\centering
\includegraphics[width=1\textwidth]{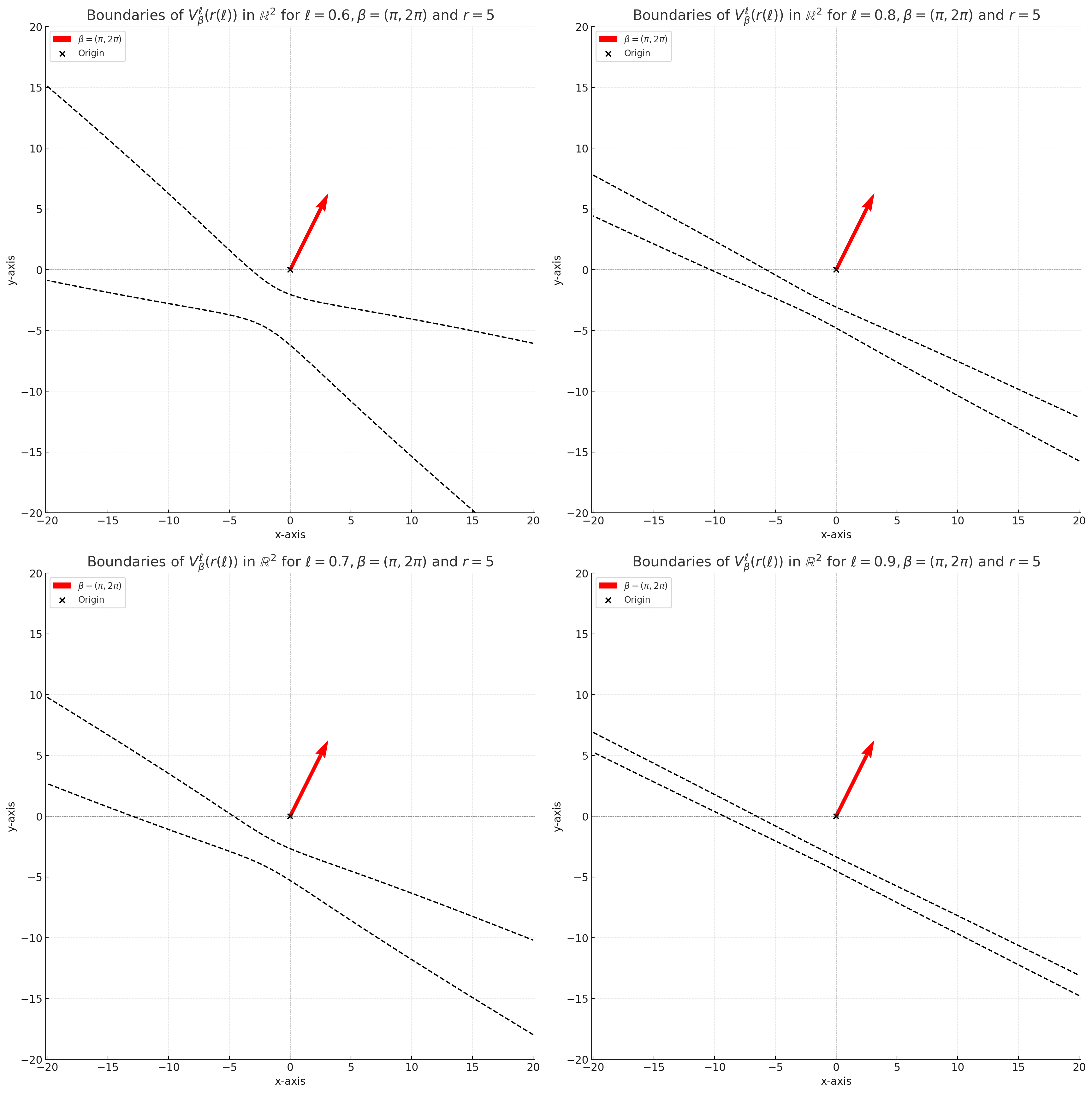}
\caption{Resonance domains for increasing values of \(\ell\).}
\label{fig:label25}
\end{figure}

\begin{figure}[htpb]
\centering
\includegraphics[width=0.4\textwidth]{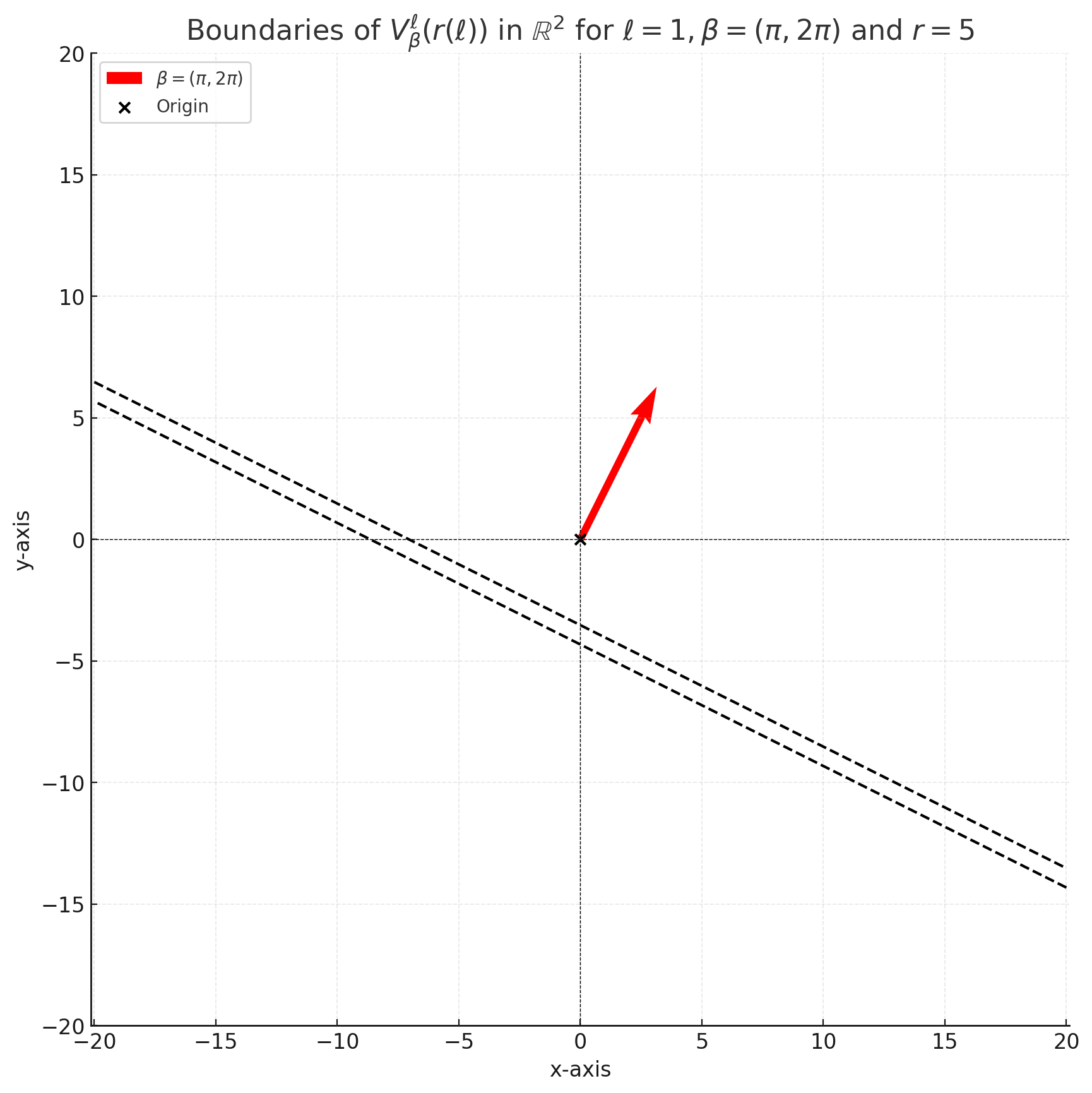}
\caption{The resonance domain for \(\ell = 1\).}
\label{fig:label}
\end{figure}

As the parameter \( \ell \) increases (e.g., \( \ell = 0.6 \to 0.9 \)), the function \( |x|^{2\ell} \) exhibits reduced sensitivity to variations in the magnitude of \( x \). This results in the following key observations:  
\begin{enumerate}
    \item \textit{A decrease in the width of the band}, since the term \( ||x|^{2\ell} - |x + \beta|^{2\ell}| \) varies more gradually.  (see Figure \ref{fig:label25}, \ref{fig:label})
    \item \textit{A contraction of the band boundaries}, leading to a more localized region. (see Figure \ref{fig:label25}, \ref{fig:label})
\end{enumerate}

\begin{figure}[htpb]
\centering
\includegraphics[width=0.8\textwidth]{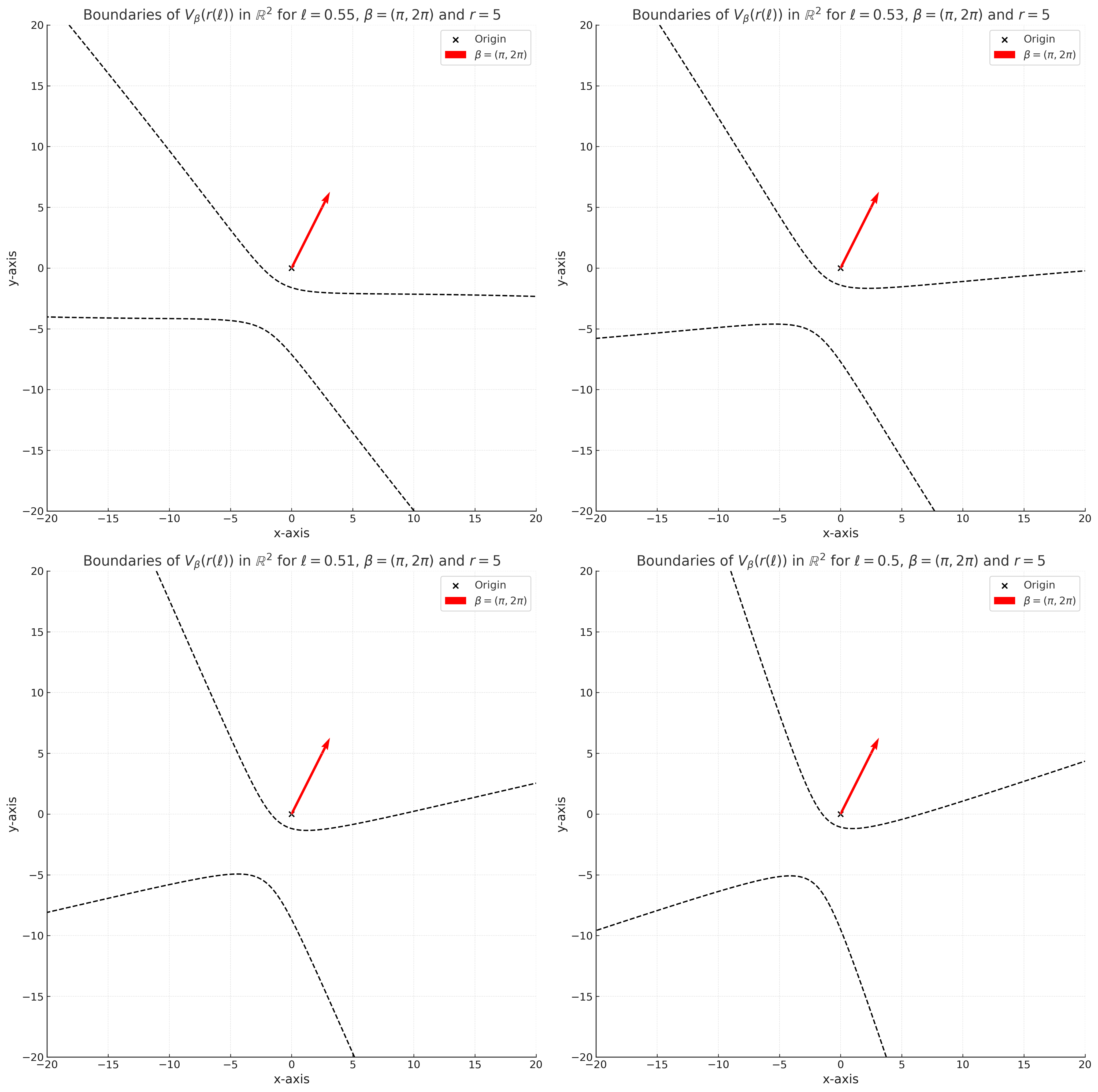}
\caption{Resonance domains for decreasing values of \(\ell\).}
\label{fig:combvisu}
\end{figure}

This effect becomes more pronounced as \( \ell \) approaches \( 1 \), at which point the boundary contours increasingly align with the differences in norms in an approximately linear manner.  On the other hand, as \( \ell \to 0.5 \), the boundaries of \( V_{\beta}^\ell(r(\ell)) \) progressively take on a parabolic shape. At \( \ell = 0.5 \), the boundaries become a perfect parabola because the norm simplifies to its linear form (\( |x| \)) (see Figure \ref{fig:combvisu}).

\begin{figure}[htpb]
\centering
\includegraphics[width=0.6\textwidth]{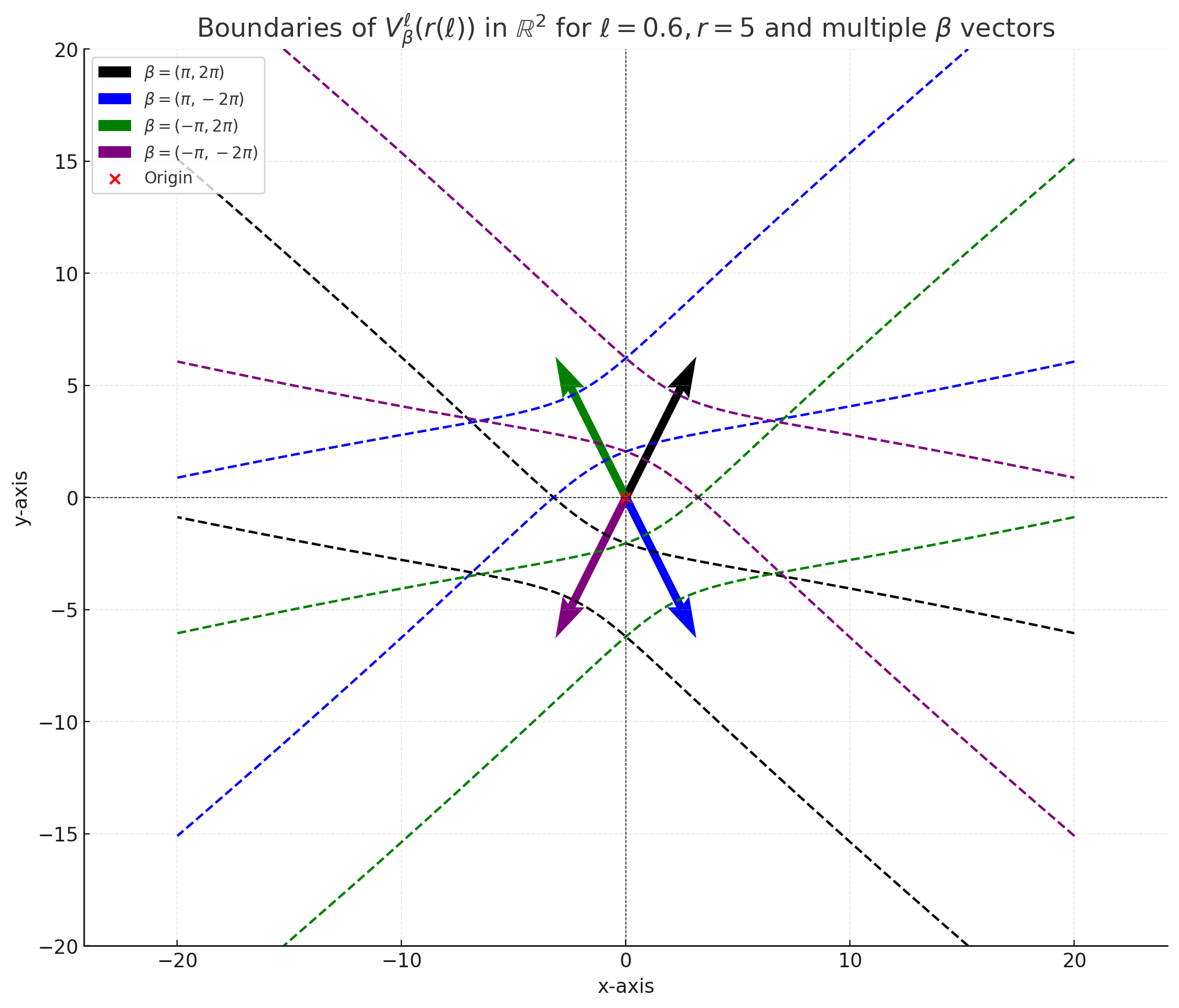}
\caption{Intersection of resonance domains corresponding to different values of \( \beta \).}
\label{inter}
\end{figure}

For different choices of \( \beta \), the resonance domains may exhibit intersections, as illustrated in Figure \ref{inter}.  

\begin{figure}[htpb]
\centering
\includegraphics[width=1\textwidth]{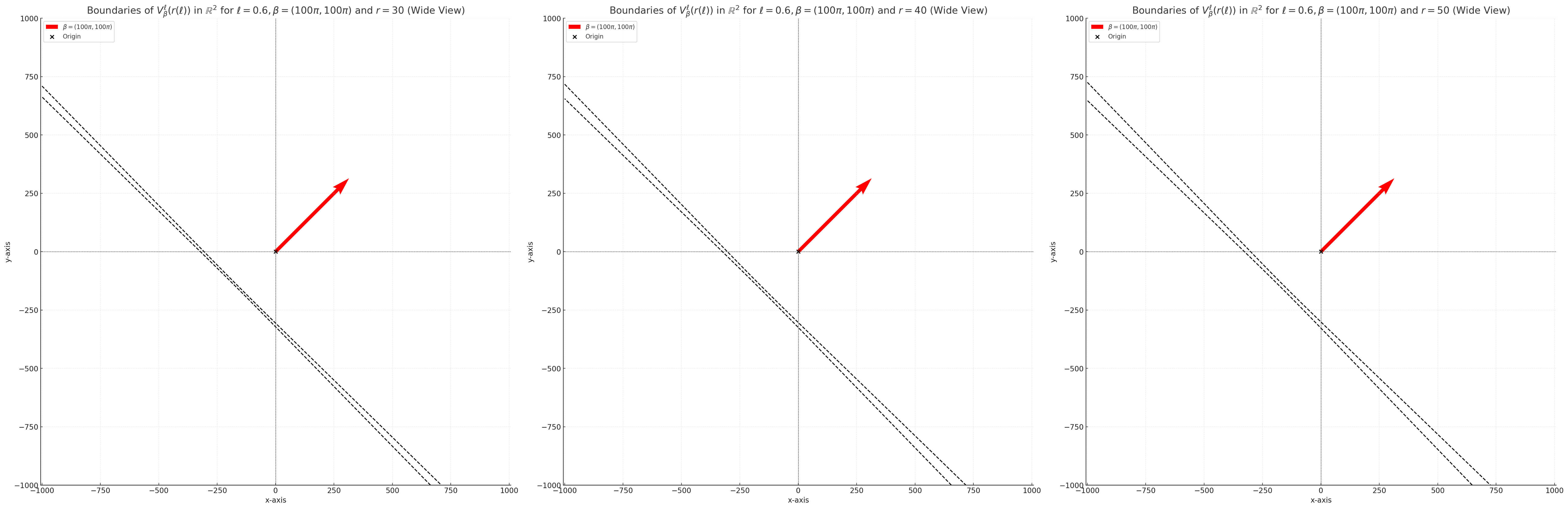}
\caption{Resonance domains for increasing values of \( r \).}
\label{rincreases}
\end{figure}

Furthermore, we analyze the effect of increasing the parameter \( r \) on the shape of these sets. As \( r \) increases, the width of the resonance domains expands, as depicted in Figure \ref{rincreases}.  

\begin{lemma}\label{mvt} For all $\beta_{1},\ldots,\beta_{k}\in \mathcal{B},$
	\begin{align}\label{Ul2}
	U(r^{\alpha_{1}(\ell)-2\ell+2},p)\subset&  
U^{\ell}(r(\ell),p)\\
\cap_{i=1}^{k}V_{\beta_{i}}^{\ell}(r^{\alpha_{k}(\ell)})\subset&\cap_{i=1}
^{k}V_{\beta_{i}}(r^{\alpha_{k}(\ell)-2\ell+2}).
	\end{align}
	\begin{proof}
		If $x \in R$,  $\left|  x\right|  \sim r$ and
		$\beta\in \mathcal{B} $  then $\left|  x+\beta\right| \sim r$. Using Mean Value Theorem we have the following equality
		\begin{equation}\label{meanvalue}
			\left|  x\right| ^{2\ell}-\left|  x+\beta\right| ^{2\ell}=\eta^{2(\ell-1)}(\left| 
			x\right| ^{2}-\left|  x+\beta\right| ^{2})
		\end{equation}
		where $\eta\sim r.$
		Hence we have
		$$
		U(r^{\alpha_{1}(\ell)-2\ell+2},p)\subset
		U^{\ell}(r(\ell),p), \qquad V_{\beta_{i}}^{\ell}(r^{\alpha_{k}(\ell)})\subset V_{\beta_{i}}(r^{\alpha_{k}(\ell)-2\ell+2}). $$
	\end{proof}
\end{lemma}
\begin{remark}
As shown in \cite{12}, using the expression \eqref{Ul2}, it can be proven that the non-resonance domain \( U^{\ell}(r(\ell),p) \) has asymptotically full measure in \(\mathbb{R}^d\). This means that the ratio \(\frac{\mu(U^{\ell}(r(\ell),p) \cap B(r))}{\mu(B(r))}\) approaches 1 as \(r \rightarrow \infty\), where \( B(r) = \{x \in \mathbb{R}^d : |x| < r\} \).
\end{remark}

\subsection{\bf Iteration Series}
In this section, we will derive the iteration series that will be used to prove the asymptotic formula.	The main condition to obtain the iteration series is the iteration condition on the vectors $\beta \in\mathcal{B}$ given by
	\begin{equation}\label{nitcon}
	|\xi_N-|\beta|^{2\ell}|>\frac{1}{2}r(\ell).
\end{equation}
This condition ensures that the "remaining terms" get small in arbitrary order after enough iterations in the asymptotic formula.  Furthermore, this condition holds  for every $\beta +\beta_{1}$ where $\beta \in  U^{\ell}(r^{\alpha_1(\ell)},p) $ is close to the eigenvalue $\xi_N$ and $\beta_{1}\in \mathcal{B}(pr^{\alpha(\ell)}) $. To see this, let's first prove the following lemma.

Let $e_i=\left( 0,\ldots,0,\frac{\pi}{a_i},0,\cdots,0\right) $ for $i=1,\ldots,d$, where $\frac{\pi}{a_i}$ lies in the $i$-th coordinate.
	By definition of the non-resonance domain, if $\beta \in U^{\ell}(r^{\alpha_1(\ell)},p)$, then $\beta \notin V^{\ell}_{e_k}(r(\ell)),$ 
	$k=1,2,\cdots,d$.
	By Lemma \ref{mvt}, in particular $\beta \notin V_{e_k}(r^{\alpha_1(\ell)-2\ell+2})$, that is,
	$$	||\beta|^{2}-|\beta+e_k|^{2}|>r^{\alpha_{1}(\ell)-2\ell +2}.$$So we have \begin{equation}\label{comp}
			|\beta^k|>\frac{1}{3}r^{\alpha_1(\ell)-2\ell+2}, \hspace{.2in} \forall
			k=1,2,\cdots,d.
		\end{equation}
		\begin{lemma}
		If $ \beta
		\in U^{\ell}(r(\ell),p) $ and $ \xi_{N}$ satisfies the inequality
		\begin{equation}\label{wellknown}
			|\xi_N-| \beta |^{2\ell}|<
			\frac{1}{2}r(\ell),
		\end{equation}
		then $ \forall\beta_{1} \in
		\mathcal{B}(pr^{\alpha(\ell)})$ the vectors $\beta + \beta_{1}$  satisfy the iteration condition
		\eqref{nitcon}.
		\begin{proof}
			Using the reverse triangle inequality, we have
	\begin{align*}
	|\xi_N-|\beta+\beta_{1}|^{2\ell}|&\geq ||\beta+\beta_{1}|^{2\ell}-|\beta|^{2\ell}|-||\beta|^{2\ell}-\xi_N|  \\
	&>r(\ell)-\frac{1}{2}r(\ell)=\frac{1}{2}r(\ell).\nonumber\qedhere
\end{align*}
		\end{proof}
	\end{lemma}  
	The following lemma plays a crucial role when we start to derive the iteration series.
	\begin{lemma} If $\beta \in\mathcal{B}$ satisfies \eqref{comp}, then 
		\begin{align} \label{qcarpim1}
			\sum_{\beta_1 \in \mathcal{B}(r^{\alpha(\ell)})}q_{\beta_1}
			v_{\beta_1}(x)v_{\beta}(x)=\sum_{\beta_1 \in
				\mathcal{B}(r^{\alpha(\ell)})}q_{\beta_1}v_{\beta+\beta_1}(x).
		\end{align}
		\begin{proof}
The proof follows a similar approach to that in \cite{18}.
		\end{proof}
	\end{lemma}
 To get the iteration series, we take inner product of both sides of $$(-\Delta_{\mathsf{NEU}})^{\ell}\chi_{N}+q(x)\chi_{N}=\xi_{N}\chi_{N} $$ with the eigenfunction $v_{\beta}$ of the Laplace operator, we obtain
	\begin{equation*}
		((-\Delta_{\mathsf{NEU}})^{\ell}\chi_{N},v_{\beta})+(q(x)\chi_{N},v_{\beta})=\xi_{N}(\chi_{N},v_{\beta}).
	\end{equation*}
As $(-\Delta_{\mathsf{NEU}})^{\ell}$ is self-adjoint and having the relation $(-\Delta_{\mathsf{NEU}})^{\ell}v_{\beta}=|\beta|^{2\ell}v_{\beta}$, we get the so-called binding formula
	\begin{align}\label{bindf1}
		(\xi_{N}-|\beta|^{2\ell} )
		(\chi_{N},v_{\beta}(x))=(\chi_{N},q(x)v_{\beta}(x) ).
	\end{align}
	We substitute the decomposition \eqref{qxorder} of $q(x)$ into the
	formula \eqref{bindf1}, then by using
	\eqref{qcarpim1},  we obtain
	\begin{align}\label{temel}
			(\xi_{N}-|\beta|^{2\ell})h(N,\beta)=& \left( \chi_{N},\sum_{\beta_1\in
				\mathcal{B}(r^{\alpha(\ell)})}q_{\beta_1}v_{\beta_{1}}v_{\beta}\right) +O(r^{-p\alpha(\ell)})\nonumber\\
	=&\sum_{\beta_1\in
			\mathcal{B}(r^{\alpha(\ell)})}q_{\beta_1}h(N,\beta+\beta_1)+O(r^{-p\alpha(\ell)}),
	\end{align}
	where we use the notation  $h(N,\beta)=(\chi_{N},v_{\beta})$. Also \eqref{bindf1}  together with \eqref{qxorder} imply
	\begin{align}\label{bngamma}
		h(N,\beta)=\frac{(\chi_{N},q(x)v_{\beta})}
		{\xi_{N}-|\beta|^{2\ell}}= \sum_{\beta_1 \in
			\mathcal{B}(r^{\alpha(\ell)})}q_{\beta_1}
		\frac{h(N,\beta+\beta_1)}
		{\xi_{N}-|\beta|^{2\ell}}+O(r^{-p\alpha(\ell)}),
	\end{align}
	for every vector $\beta\in
	{\mathcal{B}}$, satisfying the iteration condition.

	The vectors $\beta+\beta_1$ in the
	expressions \eqref{temel} satisfy the iteration condition
	\eqref{nitcon} if it is assumed that the vector $\beta\in
	U^{\ell}(r(\ell),p)$ satisfies inequality 
	\eqref{wellknown}. Under this assumption, in \eqref{bngamma}
writing $\beta+\beta_1$ and $\beta_{2}$ 	instead of $\beta$ and $\beta_{1} $, respectively, the
 expression in (\ref{bngamma}) becomes
 \begin{equation} \label{h(N,beta+beta1)}
h(N, \beta+\beta_{1})=\sum_{\beta_{2}\in \mathcal{B}(r^{\alpha(\ell)})}\frac{q_{\beta_{2}}h(N,\beta+\beta_{1}+\beta_{2})}{\xi_{N}-|\beta+\beta_1|^{2\ell}}.
 \end{equation}
 Putting the last equation into \eqref{temel}, we get the first iteration
	\begin{align*}
		(\xi_{N}-|\beta|^{2\ell} ) h(N, \beta)
		=\sum_{\beta_1,\beta_2\in \mathcal{B}(r^{\alpha(\ell)})}q_{\beta_1}
		q_{\beta_2} \frac{h(N,\beta+\beta_1+\beta_2)}
		{\xi_{N}-|\beta+\beta_1|^{2\ell}} +O(r^{-p\alpha(\ell)}).
	\end{align*}
	Isolating the terms with coefficient $ h(N,\beta) $ in the last
	equation, we obtain
	\begin{align*}
		(\xi_{N}-|\beta|^{2\ell})h(N, \beta)&= 	\sum_{{\beta_{1},\beta_{2} \in \mathcal{B}(r^{\alpha(\ell)})}
			\atop{\beta_{1}+\beta_{2}=0}}
		q_{\beta_{1}}q_{\beta_{2}}\frac{h(N, \beta)}{\xi_{N}-|\beta
			+ \beta_{1}|^{2\ell}} 
		\\ &+\sum_{{\beta_{1},\beta_{2}
				\in\mathcal{B}(r^{\alpha(\ell)})}\atop{\beta_{1}+\beta_{2}\neq 0}}q_{\beta_{1}}q_{\beta_{2}}\frac{h(N,\beta+\beta_{1}+\beta_{2}
			)}{\xi_{N}-|\beta + \beta_{1}|^{2\ell}} +O(r^{-p\alpha(\ell)})
	\end{align*}
For the next iteration, in the expression \eqref{h(N,beta+beta1)} we write  $\beta+\beta_1+\beta_{2}$ and $\beta_{3}$ 	instead of $\beta+\beta_{1}$ and $\beta_{2}$, respectively and get a similar expression for $h(N,\beta+\beta_{1}+\beta_{2}) $. Substituting this expression into the last equation and separating the terms with coefficient  $ h(N,\beta)$, we
get the second iteration. Repeating this process $ p_{1}=\left[ \frac{p+1}{3}\right]  $ times and
	isolating each time the terms with coefficient  $ h(N,\beta)$, we
	get the iteration series
	\begin{align}\label{sumsicp}
		(\xi_{N}-|\beta|^{2\ell} )h(N,\beta)
		=&\sum^{p_{1}}_{i=1}S_{i}(\xi_{N}) h(N,\beta)+
		C_{p_{1}}(\xi_N,h(N,\beta+\beta_{1}+\cdots+\beta_{p_1}))\nonumber\\
		+&O(r^{-p\alpha(\ell)})
	\end{align}
	where
	\begin{equation} \label{sumsi}
		S_{j}(\xi_{N})=\sum_{\beta_{1},\cdots,\beta_{j+1} \in
			\mathcal{B}(r^{\alpha(\ell)})} \frac{q_{\beta_{1}}\cdots q_{\beta_{j+1}}}
		{(\xi_{N}-|\beta +\beta_{1}|^{2\ell})\cdots(\xi_N-|\beta
			+\beta_{1}+\cdots+\beta_{j}|^{2\ell})},
	\end{equation}
	and
	\begin{equation}\label{dcp}
		C_{p_{1}}(\xi_{N},\zeta)\nonumber =\sum_{{\beta_{1},\cdots,\beta_{p_{1}+1}
				\in \mathcal{B}(r^{\alpha(\ell)})}\atop
			{\beta_{1}+\beta_{2}+\cdots+\beta_{p_{1}+1}\neq
				0}}\frac{q_{\beta_{1}}\cdots q_{\beta_{p_{1}+1}}\zeta}
		{(\xi_{N}-|\beta
			+\beta_{1}|^{2\ell})\cdots(\xi_{N}-|\beta +\beta_{1}
			+\cdots+\beta_{p_{1}}|^{2\ell})}.
	\end{equation}
	\begin{lemma} 
		 One has
		\begin{equation}\label{sirin}
			S_{j}(\xi_{N})=O(r^{-j\alpha_{1}(\ell)}),
		\end{equation}
		for $ j=1,2,\cdots,p_{1}$ and
		\begin{equation}\label{cip1}
			C_{p_{1}}(\xi_N,h(N,\beta+\beta_{1}+\cdots+\beta_{p_1}))=O(r^{-p_{1}\alpha_{1}(\ell)}).
		\end{equation}
		\begin{proof}
			Since the vectors $\beta_{j} \in \mathcal{B}(r^{\alpha(\ell)})$	for all $ j=1,2,\cdots,p_{1}$, we have
			$|\beta_{1}+\beta_{2}+\cdots+\beta_{j}|<p_{1} r^{\alpha(\ell)}$.
			Therefore, the iteration condition \eqref{nitcon} holds for the elements $\beta+\beta_{1}+\beta_{2}+\cdots+\beta_{j} $ for all $ j=1,2,\cdots,p_{1}$. Using \eqref{Mqx} together with Cauchy-Schwarz inequality, we see that the summation 
			 $$\sum_{\beta_{1},\cdots\beta_{j+1} \in
				\mathcal{B}(r^{\alpha(\ell)})} q_{\beta_{1}}\cdots q_{\beta_{j+1}}$$
				 is bounded. Then we conclude				\begin{align*}
	S_{j}(\xi_{N})\leq &2^{j}r^{-j\alpha_{1}(\ell)} \left(  \sum_{\beta_{1},\cdots,\beta_{j+1}\in \mathcal{B}(r^{\alpha(\ell)})}q_{\beta_{1}}\cdots q_{\beta_{i+1}}\right) \\
	=&O(r^{-j\alpha_{1}(\ell)}).
				\end{align*}
					For the latter estimate, we also need to show that $h(N,\beta+\beta_{1}+\cdots+\beta_{p_1})$ is bounded. Using Cauchy-Schwarz inequality and boundedness of the eigenfunctions $||\chi_{N}||=||v_{\beta+\beta_{1}+\cdots+\beta_{p_1}}||=O(1)$, we get \eqref{cip1}. 
		\end{proof}
	\end{lemma}
	\begin{corollary}Let $ I=\left[ |\beta|^{2\ell}-\frac{1}{2}r(\ell),|\beta|^{2\ell}+\frac{1}{2}r(\ell)\right]$. Then we have
		\begin{equation} \label{sia}
			\sum_{i=1}^{p_{1}} S_{i}(\chi)=O(r^{-\alpha_{1}(\ell)}),
			\hspace{.7cm} \forall \chi \in I.
		\end{equation}
		\begin{proof}
			To obtain \eqref{sirin}, we have only used the iteration
			condition \eqref{nitcon} and \eqref{wellknown}, that is,  $\xi_{N} \in
			I=\left[ |\beta|^{2\ell}-\frac{1}{2}r(\ell),|\beta|^{2\ell}+\frac{1}{2}r(\ell)\right] 
			$.
		\end{proof}
	\end{corollary}
	To close this section, we will prove that  for any $\beta \in\mathcal{B}$, one can find an eigenvalue $\xi_{N}$ that is close to the eigenvalue $|\beta|^{2\ell}$ and have the property that the inner product $\left| (\chi_{N},v_{\beta})\right|$ is greater than a multiple of $r^{-\frac{d-2\ell}{2}}$.
	\begin{lemma}\label{innerprodmagnitude}
	There is an eigenvalue $\xi_{N}$ satisfying the following inequalities
	\begin{equation}\label{sksk}
		|\xi_{N}-|\beta|^{2\ell}| <\frac{1}{2}r(\ell),
		\hspace{.2in} |h(N,\beta)|>c_{3}r^{-\frac{(d-2\ell)}{2}},
	\end{equation}
	where $c_3$ is a constant.
	\begin{proof} By \cite[Chapter XII.2.12, Theorem 12]{25} we know that the set of eigenfunctions $\{\chi_{N} \}_{N \in \mathbb{N}}$ of the fracrional Schrödinger operator is complete. Then if we 
		use the binding formula \eqref{bindf1}
		together with the Bessel's inequality, we have
		\begin{align*}
			\sum_{N:|\xi_{N}-|\beta|^{2\ell}|>\frac{1}{2}r(\ell)}|h(N,\beta)|^{2}&=
			\sum_{N:|\xi_{N}-|\beta|^{2\ell}|>\frac{1}{2}r(\ell)}\frac{|(\chi_{N},q(x)v_{\beta}(x))|^{2}}
			{|\xi_{N}-|\beta|^{2\ell}|^{2}}\\
			&=O(r^{-2\alpha_1(\ell)}) .
		\end{align*}
		Therefore, by Parseval's identity
		\begin{equation}\label{parseval*}
			\sum_{N:|\xi_{N}-|\beta|^{2\ell}|\leq
				\frac{1}{2}r(\ell)}|h(N,\beta)|^{2}=1-O(r^{-2\alpha_1(\ell)}).
		\end{equation}
We have the following estimate:  
\[
\left| \{\delta\in \mathbb{R}^{d}: a \sim r, \left| |\delta|^{2\ell}-a^{2\ell} \right| <1 \}\right|=O(r^{d-2\ell}).
\] Hence the number of eigenvalues of  \(\mathsf{H}_{\mathsf{NEU}}(\ell,0)\) belonging to $(a^{2\ell}-1,a^{2\ell}+1)$ has order $O(r^{d-2\ell})$. By \cite[Chapter V, Section 4.3]{21} we know that the \(N\)-th eigenvalue of \(\mathsf{H}_{\mathsf{NEU}}(\ell,q)\) is within an \(M\)-neighborhood of the \(N\)-th eigenvalue of \(\mathsf{H}_{\mathsf{NEU}}(\ell,0)\). So the number of eigenvalues \( \xi_{N} \) of the operator \( \mathsf{H}_{\mathsf{NEU}}(\ell,q) \) contained in the interval $\left[ |\beta|^{2\ell}-M,|\beta|^{2\ell}+M\right] $ is of order \( O(r^{d-2\ell}) \).  As \[
\left[ |\beta|^{2\ell}-M,|\beta|^{2\ell}+M\right] \subseteq \left[ |\beta|^{2\ell}-\frac{1}{2}r(\ell),|\beta|^{2\ell}+\frac{1}{2}r(\ell)\right],
\]  
		there exists $\xi_ N \in I$  such that $$|h(N,\beta)|^{2}> \frac{1-O(r^{{-2\alpha_{1}(\ell)}})}{O(r^{d-2\ell})}=c_{3}r^{-(d-2\ell)}.\eqno\qedhere$$
	\end{proof}
\end{lemma}
\subsection{\bf Asymptotic Formula}\label{asymp}

In this section, we will obtain the asymptotic formula for the eigenvalues of the operator $\mathsf{H}_{\mathsf{NEU}}.$ To do this, we will use the results that we obtain on the terms of iteration series in the previous section. Also we will introduce the following sequence.	For $k=2,3,\cdots,p-c$
	and $c=\left[ \frac{d-1}{2\alpha(\ell)}\right] +1 $, we define 
	\begin{align*}
F_{0}=0,& \hspace{.2cm} F_{1}=\sum_{\beta_{1} \in
	\mathcal{B}(r^{\alpha(\ell)})} \frac{|q_{\beta_{1}}|^{2}}{
	|\beta|^{2\ell}-|\beta+\beta_{1}|^{2\ell}},\\ &F_{k}=\sum_{i=1}^{k}S_{i}(|\beta|^{2\ell}+F_{k-1}).
	\end{align*}
	\begin{lemma}\label{fjest} For all $j=1,\ldots, p_1$
		\begin{equation}
			F_{j}=O(r^{-\alpha_{1}(\ell)}).
		\end{equation}
		\begin{proof}
			We prove by induction on $k$. For the initial case $F_0=0.$ Suppose $F_{j-1} =
			O(r^{-\alpha_{1}(\ell)}) $, then $|\beta|^{2\ell}+F_{j-1} \in I$ and by
			(\ref{sia}) we have $
			F_j=\sum_{i=1}^{j}S_{i}(|\beta|^{2\ell}+F_{j-1})=O(r^{-\alpha_{1}(\ell)})$.
		\end{proof}
	\end{lemma}
	\begin{theorem}
		Let $\beta \in U^{\ell}(r^{\alpha(\ell)},p) $  and
		$|\beta|\sim r$.  Then there is an eigenvalue $ \xi_{N}$ of
		the operator $\mathsf{H}_{\mathsf{NEU}}(\ell,q) $ satisfying the formula
		\begin{equation}\label{fk1}
			\xi_{N}=|\beta|^{2\ell} + F_{k-1} +O(r^{-k\alpha_1(\ell)}),
		\end{equation}
		for all $k=1,2,\cdots,p-c $.
		\begin{proof}
			We prove (\ref{fk1}) by induction on $k$. For $k=1$, we substitute (\ref{sirin}) into the iteration series (\ref{sumsicp}) and we get
			\begin{align}\label{last}
				(\xi_{N}-|\beta|^{2\ell} ) h(N, \beta)  =
				O(r^{-\alpha_{1}(\ell)})h(N, \beta)+
				O(r^{-p\alpha(\ell)}).
			\end{align}
By Lemma \ref{innerprodmagnitude}, there exists an eigenvalue \(\xi_N\) satisfying the following properties:  
\begin{equation} \label{indk=1}
    \left| \xi_N - \|\beta\|^{2\ell} \right| < \frac{1}{2} r(\ell),  
    \hspace{.2in} |  h(N, \beta)| > C r^{-\frac{(d-2\ell)}{2}},
\end{equation}
where \(C\) is a constant. Dividing both sides of equation \eqref{last} by \((\chi_N, v_{\beta})\) for this eigenvalue \(\xi_N\) satisfying \eqref{indk=1}, we obtain  
			\begin{equation*}
 	\xi_{N}=|\beta|^{2\ell} +O(r^{-\alpha_{1}(\ell)}).
			\end{equation*}
			  Assume that (\ref{fk1}) holds for $k=j$,
			that is;
			\begin{equation}\label{assumej}
				\xi_{N}=|\beta|^{2\ell}+F_{j-1} +O(r^{-j\alpha_{1}(\ell)}).
			\end{equation}
			To prove  for $ k=j+1 $, we substitute the expression for $\xi_{N}$
			in \eqref{assumej} into the term $S_{i}(\xi_{N})$ of the
			formula \eqref{sumsicp}. Then we get
			\begin{equation}\label{*1}
				\footnotesize
			(	(\xi_{N}-|\beta|^{2\ell} ) h(N, \beta)
				= \sum^{p_{1}}_{i=1}S_{i}(|\beta|^{2\ell}+F_{j-1}
				+O(r^{-j\alpha_{1}(\ell)}))) h(N, \beta)+ C_{p_{1}}+
				O(r^{-p\alpha(\ell)}) ),
			\end{equation}
			dividing  both sides of the  equation (\ref{*1}) by $h(N, \beta)$
			and then  using (\ref{sksk}) and (\ref{cip1}), we get
			\begin{equation}\label{sumN}
				\xi_{N}=|\beta|^{2\ell} +\sum_{i=1}^{p_{1}}
				S_{i}(|\beta|^{2\ell}+F_{j-1}+O(r^{-j\alpha_{1}(\ell)}))+O(r^{-(p-c)\alpha(\ell)})
			\end{equation}
			Adding and subtracting the term
			$F_j=\sum_{i=1}^{j}S_{i}(|\beta|^{2\ell}+F_{j-1})$ in (\ref{sumN}), we have
			{\footnotesize\begin{align}\label{sumnbr}
					\xi_{N}&=|\beta|^{2\ell} +\left[ \sum_{i=1}^{j}
					S_{i}(|\beta|^{2\ell}+F_{j-1}+O(r^{-j\alpha_{1}(\ell)}))-
					S_{i}(|\beta|^{2\ell}+F_{j-1}) \right]  \nonumber\\
					&+\sum_{i=j+1}^{p_{1}}
					S_{i}(|\beta|^{2\ell}+F_{j-1}+O(r^{-j\alpha_{1}(\ell)}))+\sum_{i=1}^{j}
					S_{i}(|\beta|^{2\ell} +F_{j-1})+ O(r^{-(p-c)\alpha(\ell)}).
			\end{align}}
			Notice that, (\ref{sia}) implies
			$$\sum_{i=j+1}^{p_{1}}
			S_{i}(|\beta|^{2\ell}+F_{j-1}+O(r^{-j\alpha_{1}(\ell)}))=O(r^{-(j+1)\alpha_{1}(\ell)}).$$
			So, we need only to show that the expression in the square brackets
			in (\ref{sumnbr})  is equal to $
			O(r^{-(j+1)\alpha_{1}(\ell)}) $. Using Lemma \ref{fjest}, we see that $|\beta|^{2\ell}
			+F_{j-1}+O(r^{-j\alpha_1(\ell)}),|\beta|^{2\ell}
			+F_{j-1}\in I$. So  the iteration condition gives
            $$ 		\bigg||\beta|^{2\ell}
			+F_{j-1}+O(r^{-j\alpha_1(\ell)})-
	|\beta+\beta_{1}+\cdots+\beta_{i}|^{2\ell}\bigg|>\frac{1}{3}r(\ell)$$ and $$ \bigg| |\beta|^{2\ell} +F_{j-1}-
			|\beta+\beta_{1}+\cdots+\beta_{i}|^{2\ell}\bigg| >\frac{1}{3}r(\ell),
			$$ for all $i=1,2,\cdots,p_{1}.$
			Using the last two inequalities, \eqref{Mqx} and by direct calculations, the expression in the square brackets is $O(r^{-(j+1)\alpha_{1}(\ell)})$:  
			{\scriptsize\begin{align*}
					&\frac{1}{ \prod_{j=1}^{s}\left( |\beta|^{2\ell}+F_{j-1}+O(r^{-j\alpha_{1}})-\left|\beta+\sum_{i=1}^{j}\beta_{i} \right|^{2\ell} \right) }-\frac{1}{\prod_{j=1}^{s}\left( |\beta|^{2\ell}+F_{j-1}-\left|\beta+\sum_{i=1}^{j}\beta_{i} \right|^{2\ell} \right) }\\
					=&\frac{1}{\prod_{j=1}^{s}\left( |\beta|^{2\ell}+F_{j-1}-\left|\beta+\sum_{i=1}^{j}\beta_{i} \right|^{2\ell}\right)  } \left(\frac{1}{1-O(r^{-(j+1)\alpha_{1}(\ell)})} \right)\\=&O(r^{-(j+1)\alpha_{1}(\ell)}).\qedhere
			\end{align*}}
		\end{proof}
	\end{theorem}

\end{document}